%% file: main.tex
\documentclass{ifacconf}
\input{header}

\usepackage{natbib} 

\def\pars#1{{\color{black}#1}} 
\def\pd#1{{\color{black}#1}} 
\def\ar#1{{\color{black}#1}}

\begin{document}
\begin{frontmatter}

\title{A General Framework for Distributed Partitioned Optimization} 

\thanks[footnoteinfo]{The research was supported by Russian Science Foundation (project No. 21-71-30005).}
\thanks[equal_contribution]{S. Chezhegov and A. Novitskii contributed equally.}

\author[ISPRAS]{Savelii Chezhegov}
\author[ISPRAS]{Anton Novitskii} 
\author[MIPT]{Alexander Rogozin}
\author[Skoltech]{Sergei Parsegov}
\author[WIAS]{Pavel Dvurechensky}
\author[MIPT,ISPRAS]{Alexander Gasnikov}

\address[ISPRAS]{ISP RAS Research Center for Trusted Artificial Intelligence, Moscow, Russia}
\address[MIPT]{Moscow Institute of Physics and Technology, Dolgoprudny, Russia}
\address[Skoltech]{Skolkovo Institute of Science and Technology, Moscow, Russia}
\address[WIAS]{Weierstrass Institute for Applied Analysis and Stochastics, Berlin, Germany}

\begin{abstract} 
Distributed optimization is widely used in large scale and privacy preserving machine learning and various distributed control and sensing systems. It is assumed that every agent in the network possesses a local objective function, and the nodes interact via a communication network. In the standard scenario, which is mostly studied in the literature, the local functions are dependent on a common set of variables, and, therefore, have to send the whole variable set at each communication round. In this work, we study a different problem statement, where each of the local functions held by the nodes depends only on some subset of the variables. 
Given a network, 
we build a general algorithm-independent framework for distributed partitioned optimization that allows to construct algorithms with reduced communication load using a generalization of Laplacian matrix. Moreover, our framework allows to obtain algorithms with non-asymptotic convergence rates with explicit dependence on the parameters of the network, including accelerated and optimal first-order methods. We illustrate the efficacy of our approach on a synthetic example.
\end{abstract}

\begin{keyword}
Large scale optimization problems, Convex optimization, Optimization and control of large-scale network systems, Decentralized and distributed control, Multiagent systems
\end{keyword}

\end{frontmatter}







\input{sections/introduction}
\input{sections/problem}

\input{sections/ring_of_cliques}
\input{sections/conclusion}

\bibliography{references}

\newpage
\onecolumn
\appendix
\input{sections/appendix}

\end{document}

%% file: header.tex
\usepackage{url}            
\usepackage{booktabs}       
\usepackage{amsfonts}       
\usepackage{nicefrac}       
\usepackage{microtype}      
\usepackage{lipsum}
\usepackage{mathtools}
\usepackage{subcaption}

\usepackage{amsmath,amssymb}
\usepackage{pifont}
\usepackage{cases}
\usepackage{stackengine}    
\usepackage{wrapfig}

\usepackage{xcolor}

\newcommand{\eps}{\varepsilon}

\newcommand{\ol}{\overline}
\newcommand{\one}{\mathbf{1}}

\renewcommand{\hat}{\widehat}
\newcommand{\cd}{{\cdots}}
\newcommand{\vd}{{\vdots}}

\DeclareMathOperator{\diag}{diag}

\newcommand{\R}{\mathbb{R}}

\newcommand{\mA}{{\bf A}}
\newcommand{\mB}{{\bf B}}

\newcommand{\mI}{{\bf I}}

\newcommand{\mL}{{\bf L}}

\newcommand{\mS}{{\bf S}}

\newcommand{\mW}{{\bf W}}

\newcommand{\cE}{{\mathcal{E}}}

\newcommand{\cG}{{\mathcal{G}}}

\newcommand{\cN}{{\mathcal{N}}}

\newcommand{\cV}{{\mathcal{V}}}

\newcommand{\be}{{\bf e}}

\newcommand{\bx}{{\bf x}}

\newcommand{\bLambda}{{\bf \Lambda}}

\newcommand{\ds}{\displaystyle}

\newcommand{\cbraces}[1]{\left( #1 \right)}
\newcommand{\sbraces}[1]{\left[ #1 \right]}
\newcommand{\braces}[1]{\left\{ #1 \right\}}

\newcommand{\ar}[1]{\textcolor{brown}{#1}}
\newcommand{\pd}[1]{\textcolor{magenta}{#1}} 

\definecolor{green}{rgb}{0.3, 0.7, 0.2}
\newcommand{\blue}[1]{{\color{blue}{#1}}}
\newcommand{\green}[1]{{\color{green}{#1}}}

\def\R{\mathbb R}

\def\one{{\mathbf 1}}

\def\<#1,#2>{\langle #1,#2\rangle}

%% file: sections/introduction.tex
\section{Introduction}
Distributed algorithms is \ar{a classical} \cite{bor82,tsi84,deg74}, yet actively developing research area with many applications including robotics, resource allocation, power system control, control of drone or satellite networks, distributed statistical inference and optimal transport, multiagent reinforcement learning \cite{xia06,rab04,ram2009distributed,kra13,uribe2018distributed,kroshnin2019complexity,ivanova2020composite}.
Recent surge of interest to such problems in optimization and machine learning is motivated by large-scale learning problems with privacy constraints and other challenges such as data being produced or stored distributedly \cite{bot10,boy11,ned15}. 
An important part of this research studies distributed optimization algorithms over arbitrary networks of computing agents, e.g. sensors or computers, which is represented by a connected graph in which two agents can communicate with each other if there is an edge between them. This imposes communication constraints and the goal of the whole system~\cite{ned09} is to cooperatively minimize a global objective using only local communications between agents, each of which has access only to a local piece of the global objective. 

In this paper, we further exploit additional structure in such problems and consider distributed  partitioned  optimization problems also known as optimization with overlapping variables and distributed optimization of partially separable objective functions. Such problems arise in many modern big-data applications, e.g., distributed matrix completion, distributed estimation in power networks, network utility maximization, distributed resource allocation, cooperative localization in wireless networks, building maps by robotic networks \cite{erseghe2012distributed,kekatos2012distributed,carli2013distributed,notarnicola2017distributed,cannelli2020asynchronous}. As in the standard formulations, the goal is to minimize a sum of $m$ functions $\tilde{f}_i(x)$, $i=1,...,m$ with each $\tilde{f}_i(x)$ stored at a node of the computational network. Yet, unlike standard distributed problems, the space of decision variables is divided into $n$ blocks, and each of $\tilde{f}_i(x)$ may depend only on a, possibly small, subset of blocks. Such sparse structure leads to inefficiency of the standard approaches \cite{scaman2017optimal,kovalev2020optimal} since they require each node to store and send the whole vector of variables instead of storing and exchanging with other nodes a small vector of variables which influences its local objective $\tilde{f}_i(x)$. This leads to inefficient usage of computational and communication resources.

Theory of algorithms exploiting such additional structure seems to be underdeveloped in the literature. 
\cite{necoara2016parallel} propose a parallel version of a randomized (block) coordinate descent method for minimizing the sum of a partially separable smooth convex function and a fully separable non-smooth convex function. Moreover, they explain how to implement their algorithms in a distributed setup and obtain convergence rate guarantees. \cite{cannelli2020asynchronous} consider convex and nonconvex constrained optimization with a partially separable objective function and propose an asynchronous algorithm with rate guarantees for this
class of problems. Finally, \cite{notarnicola2017distributed} propose asynchronous dual decomposition algorithm for such problems and prove its asymptotic convergence. 

Despite very advanced results and techniques, these works have two limitations. Firstly, their distributed algorithms assume that the number of functions $m$ is the same as the number of variables blocks $n$ and each node $i$ stores not only the objective $\tilde{f}_i$, but also the $i$-th block of variables. Secondly, and more importantly, they assume that the computational graph is aligned with the dependence of $\tilde{f}_i$'s on the blocks of variables. The latter means that if $\tilde{f}_i$ depends on the block variable $x^\ell$, then the nodes $i$ and $\ell$ are connected by an edge of the computational network. In this paper, we do not make such assumptions and consider a more general setting. Moreover, we propose a general algorithm-independent framework that allows to reformulate distributed  partitioned  optimization problem using in a way suitable for application of many distributed optimization algorithms, see, e.g., \cite{yang2019survey,gorbunov2020recent,dvinskikh2021decentralized}. Our approach makes it possible to go beyond optimization problems and apply it to distributed methods for saddle-point problems \cite{rogozin2021decentralized} and variational inequalities \cite{kovalev2022optimal}. At the core of our framework lie communication matrices, e.g., Laplacian of the computational graph, which are widely used in distributed optimization \cite{gorbunov2020recent}. Our approach allows to flexibly choose computational subgraph for each block of variables and the corresponding mixing matrices in order to make the storage, computational, and communication complexity smaller.
This, in particular, allows obtaining algorithms with non-asymptotic convergence rates (unlike \cite{notarnicola2017distributed}) with explicit dependence on the parameters of the computational network, including accelerated and optimal algorithms (unlike \cite{cannelli2020asynchronous,necoara2016parallel}).

We illustrate the effectiveness of our reformulation-based framework by considering graphs of a certain structure. These graphs have a two-layer hierarchy: the first layer is represented by groups of nodes that communicate on their local variable blocks, while the upper level reflects the communications between the groups. Within the inter-group information exchange, the nodes from each group share the variables according to the links between the groups. The reasons to consider such topologies are as follows: first, such structures are scalable in terms of the number of groups and the number of agents within each group. Second, these graphs admit closed-form calculation of the Laplacian spectra, which influence the convergence rates of distributed algorithms. Moreover, such hierarchical graphs mimic the nature of the distributed estimators that consider local variables as private information and exchange the shared variables with certain neighbors only (e.g. in distributed power system state estimation, see \cite{kekatos2012distributed,notarnicola2017distributed}.
Such graphs allow us to study the asymptotics of the condition number of the Laplacian matrix for both approaches: with a common state vector and with blocks of variables, and show that our approach leads to better convergence rates of distributed algorithms.

\ar{
This paper is organized as follows. In Section \ref{sec:distr_part_optimization}, we describe our framework for partitioned optimization problems. Namely, we describe the generalization of Laplacian matrix in Section \ref{subsec:complex_lapl_framework}, give an example of building such a matrix in Section \ref{subsec:complex_lapl_example} and analyze its spectral characteristics in Section \ref{subsec:complex_lapl_spectrum}. After that, we illustrate how our approach works on a synthetic network example in Section \ref{sec:ring_of_cliques}. There we consider a hierarchical network that consists of $n$ cliques of size $k$ connected by a ring graph and show that using our approach decreases the condition number of the communication matrix by $\Theta(n^2 k)$ times.
}

\subsection{Notation}

Throughout this paper, $\mL(\cG)$ denotes the Laplacian matrix of graph $\cG = (\cV, \cE)$: 
\begin{align*}
    [\mL(\cG)]_{ij} = 
    \begin{cases}
    \text{deg}(i), &\text{if } i = j, \\
    -1, &\text{if } (i, j)\in\cE, \\
    0, &\text{else},
    \end{cases}
\end{align*}
where $\text{deg}(i)$ denotes the number of nodes adjacent to node $i$. We also let $\mI_p$ be the identity matrix of size $p\times p$, $\one_p$ be \pd{the all-ones} vector of length $p$ and $\bf{0}_p$ be a vector consisting of $p$ zeros. We denote $\be_p^{(q)} = (0, \ldots, 1, \ldots, 0)^\top$ the $q$-th coordinate vector in $\R^p$. \ar{After that, $\lambda_{\max}(\cdot)$ and $\lambda_{\min}^+(\cdot)$ are the largest and smallest positive eigenvalues of matrix, respectively, and $\chi(\cdot) = \lambda_{\max}(\cdot) / \lambda_{\min}^+(\cdot)$ is the condition number}. We also denote $A\otimes B$ a Kronecker product of matrices $A$ and $B$. Finally, $\bLambda(\cdot)$ denotes the set of \textit{unique} eigenvalues of a matrix.

%% file: sections/problem.tex
\section{Distributed partitioned optimization}\label{sec:distr_part_optimization}

\subsection{\pd{The proposed framework}}
\label{subsec:complex_lapl_framework}

Consider the following distributed partitioned optimization problem:
\begin{align}\label{eq:problem_initial}
	\min_{x\in \R^n} f(x),~~f(x) := \sum_{i=1}^m  \tilde f_i\cbraces{x^{[\cN_i]}},
\end{align}
where $\cN_i \subseteq \braces{1, \ldots, n}$, i.e., each function $\tilde f_i$ depends on a subset $x^{[\cN_i]}$ of variables\footnote{For simplicity, we consider variables $x^\ell\in \R$, but everything can be straightforwardly generalized for the case when $x^\ell \in \R^{n_\ell}$ are blocks of variables which do not intersect and the sum of all their dimensions is equal to $n$.}
$x^\ell\in \R$, $\ell = 1,\ldots,n$ and these subsets  may be of different size and even overlap. 
Further, we assume that there is a computational network represented by a connected graph $\cG = (\cV, \cE)$, where $\cV = \braces{1, \ldots, m}$ is the set of nodes and $\cE$ is the set of edges connecting the elements of $\cV$, and that each $\tilde f_i$ is locally held by a separate computational node of $\cG$. The goal of the network is to cooperatively solve problem \eqref{eq:problem_initial} under communication constraints: two nodes may exchange information if and only if there is an edge in $\cE$ connecting these nodes. Unlike previous works, it is allowed that two functions $\tilde f_i$ and $\tilde f_j$ depend on the same variable $x^\ell$, but nodes $i$ and $j$ are not connected by an edge in $\cG$.

To exploit the partitioned structure of the problem, for every variable $x^\ell$, we define a set of nodes that hold functions dependent on $x^\ell$: $\cV^\ell\subseteq\cV$ (i.e. $\cV^\ell = \braces{i:~ \ell\in\cN_i}$) and consider an undirected and connected communication subnetwork $\tilde\cG^\ell = (\cV^\ell, \cE^\ell)$ with $\cE^\ell\subseteq\cE$. By construction, an edge $(i, j)$ lies in $\cE^\ell$ if $\tilde f_i$ and $\tilde f_j$ depend on $x^\ell$ and $(i, j) \in \cE$. Note that it is possible that functions $\tilde f_i, \tilde f_j$ depend on $x^\ell$, but $(i, j)\notin\cE$. In this case $(i, j)\notin\cE^\ell$, but still the nodes $i$ and $j$ are connected by some path in $\tilde\cG^\ell$.

The standard approaches \cite{gorbunov2020recent} to solve distributed optimization problems require each node to store a local approximation of the \textit{whole} vector $x \in \R^n$ and communicate it to the neighbors. For each node $i$, we define a vector $x_i\in \R^n$ and introduce the stacked vector $\bx = [x_1^\top \ldots x_m^\top]^\top \in\R^{mn}$. To reduce the storage requirements and the amount of communicated information, we assume that each node $i$ holds an approximation $x_i^{\ell}$ \textit{only} of the variables $x^{\ell}$ such that $\ell \in \cN_i$, i.e. $\tilde{f}_i$ depends on $x^{\ell}$. To obtain a problem equivalent to \eqref{eq:problem_initial} we impose consensus constraints $x_{j_1}^\ell = x_{j_2}^\ell = \ldots = x_{j_{|\cV^\ell|}}^\ell$, where $j_1, \ldots, j_{\cV^\ell}$ are the nodes of $\cV^\ell$. Since all the graphs $\tilde\cG^\ell$ are connected and have vertices $\cV^\ell$, we can equivalently rewrite \eqref{eq:problem_initial} using the approximations $x_i^{\ell}$ as 
\begin{align}\label{eq:problem_flattened}
    \min_{\bx\in\R^{mn}}~ &F(\bx) = \sum_{i=1}^m \tilde f_i(x_i^{[\cN_i]}) \\
    \text{s.t.}~ &x_i^\ell = x_j^\ell~~ \forall (i, j)\in\cE^\ell,~ \ell = 1,\ldots, n. \nonumber
\end{align}

The next reformulation step is based on stating the constraints of this problem as a system of linear equations using a communication matrix, which requires some notation. 
We let $\bx^\ell = [x_1^\ell \ldots x_m^\ell]^\top\in\R^m$ be the vector of approximations to the variable $x^\ell$.
We further introduce graphs $\cG^\ell = (\cV, \cE^\ell)$ for $\ell = 1, \ldots, m$, which are graphs $\tilde\cG^\ell$ augmented with isolated vertices (if needed).
The communication matrix $\mW$ associated with the set of networks $\braces{\cG^\ell}_{\ell=1}^n$ is defined as\footnote{
Instead of $\mL(\cG^\ell)$ we can take a doubly stochastic mixing matrix.}
\begin{align}\label{eq:complex_lapl}
    \mW = \sum_{\ell=1}^n \mL(\cG^\ell) \otimes \be_n^{(\ell)} \be_n^{(\ell)\top}.
\end{align}
According to the definition of $\bx$, we have $(\mL(\cG^\ell)\otimes {\be_n^{(\ell)}} \be_n^{(\ell)\top})\bx = (\mL(\cG^\ell) \bx^\ell)\otimes {\be_n^{(\ell)}}$ and $\mW\bx = \sum_{\ell=1}^m (\mL(\cG^\ell) \bx^\ell)\otimes {\be_n^{(\ell)}}$.
Thus, we conclude that the linear constraint $\mW\bx = 0$ is equivalent to $\braces{\mL(\cG^\ell) \bx^\ell = 0}_{\ell=1}^m$, which, in turn, by the definition of the graph Laplacian, are equivalent to the constraints in \eqref{eq:problem_flattened}. Finally, we introduce $f_i(x) = \tilde f_i(x^{[\cN_i]})$ for $i = 1, \ldots, m$. Functions $f_i$ depend only on the variable subset $\cN_i$, but formally take the whole variable vector as their argument. Combining everything together, we obtain the following equivalent reformulation of \eqref{eq:problem_initial} and \eqref{eq:problem_flattened}
\begin{align}\label{eq:problem_affine_contraints}
    \min_{x\in\R^n}~ F(\bx) := \sum_{i=1}^m f_i(x_i)~ \text{ s.t. } \mW\bx = 0.
\end{align}
Thus, our framework results in this reformulation which is standard for distributed optimization and allows to apply a long list of algorithms to solve \eqref{eq:problem_initial} by solving \eqref{eq:problem_affine_contraints}. Any distributed consensus-based optimization algorithms \cite{gorbunov2020recent}, including the state-of-the-art primal (OPAPC \cite{kovalev2020optimal}) and dual (MSDA \cite{scaman2017optimal}) methods can be applied to our problem reformulation. The communication complexity of these algorithms explicitly depends on the parameters of the network through the condition number $\chi(\mW)$ of $\mW$. Thus, in what follows we focus on studying the spectrum of $\mW$.

\begin{remark}\label{rem:lapl_}
When all $f_i$ depend on the common set of variables, i.e. $\cN_i = \braces{1, \ldots, n}$ for $i = 1, \ldots, m$, we have $\mL(\cG^\ell) = W$ for $\ell = 1, \ldots, n$. In this case $\mW = W\otimes \sum_{\ell=1}^n {\be_n^{(\ell)}} \be_n^{(\ell)\top} = W\otimes\mI_n$, which is the standard communication matrix used in distributed optimization.
\end{remark}

\subsection{Example}\label{subsec:complex_lapl_example}

Consider a graph $\cG = \cbraces{\cV, \cE}$, where $\cV = \braces{1, 2, 3}, \cE = \braces{(1, 2), (1, 3)}$, i.e. $\cG$ and let us have \pars{two variables $x^1$, $x^2$}. Let \pars{$\tilde f_1 = \tilde f_1(x^1, x^2)$ }, \pars{$\tilde f_2 = \tilde f_2(x^1)$} and \pars{$\tilde f_3 = \tilde f_3(x^2)$}. We build corresponding Laplacians for the variables $x^1$ and $x^2$ as follows.

\begin{figure}[h!]
     \centering
     \begin{subfigure}[t]{0.3\linewidth}
         \centering
         \includegraphics[width=0.6\linewidth]{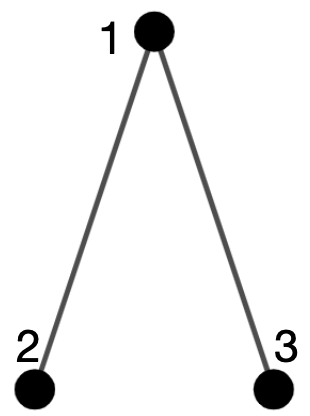}
         \caption{ Original graph $\cG$}
         \label{fig:fullgraph}
     \end{subfigure}
     \hfill
     \begin{subfigure}[t]{0.3\textwidth}
         \centering
         \includegraphics[width=0.7\linewidth]{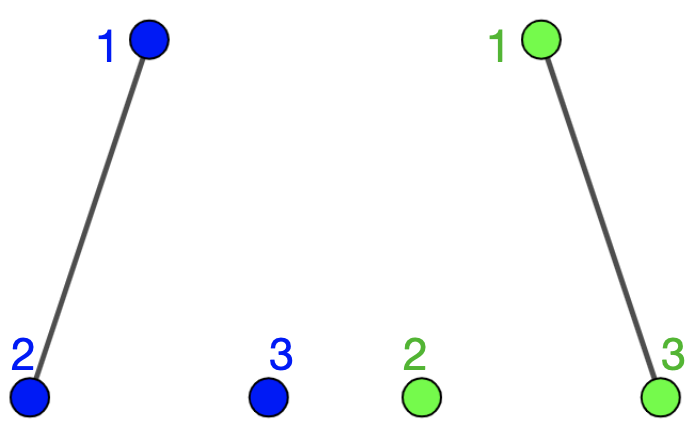}
         \caption{Subgraphs \blue{$\cG^1$} and \green{$\cG^2$}}
         \label{fig:subgraph}
     \end{subfigure}
        \caption{Graph $\cG$ \pars{of the computational network with three nodes and subgraphs $\cG^1$ (in blue), $\cG^2$ (in green)}}
        \label{fig:three graphs}
\end{figure}
\begin{align*}
\blue{
\mL(\cG^1) = 
\begin{bmatrix*}[r]
    &1&\text{-}1&0 \\
    &\text{-}1&1&0 \\
    &0&0&0 \\
\end{bmatrix*}
},\quad
\green{
\mL(\cG^2) = 
\begin{bmatrix*}[r]
    &1&0&\text{-}1 \\
    &0&0&0 \\
    &\text{-}1&0&1 \\
\end{bmatrix*}
}
\end{align*}
According to \eqref{eq:complex_lapl} $\mW$ \pars{has the following form}:
\begin{align*}
\mW &= \blue{\mL(\cG^1)}\otimes \be_1\be_1^\top + \green{\mL(\cG^2)}\otimes \be_2 \be_2^\top \\
&= \blue{
\begin{bmatrix*}[r]
    &1&\text{-}1&0 \\
    &\text{-}1&1&0 \\
    &0&0&0 \\
\end{bmatrix*}
}\otimes
\begin{bmatrix*}[r]
    &1&0 \\
    &0&0
\end{bmatrix*} + 
\green{
\begin{bmatrix*}[r]
    &1&0&\text{-}1 \\
    &0&0&0 \\
    &\text{-}1&0&1 \\
\end{bmatrix*}
}\otimes
\begin{bmatrix*}[r]
    &0&0 \\
    &0&1
\end{bmatrix*} \\
&= 
\begin{bmatrix*}[r]
    &\blue{1} &0 &\blue{\text{-}1} &0 &\blue{0} &0 \\
    &0 &\green{1} &0 &\green{0} &0 &\green{\text{-}1} \\
    &\blue{\text{-}1} &0 &\blue{1} &0 &\blue{0} &0 \\
    &0 &\green{0} &0 &\green{0} &0 &\green{0} \\
    &\blue{0}&0&\blue{0}&0&\blue{0}&0 \\
    &0 &\green{\text{-}1} &0 &\green{0} &0 &\green{1} \\
\end{bmatrix*}
\end{align*}

On the other hand, the Laplacian of the original network $\cG$ (not taking into account the different variables) writes as 
\begin{eqnarray*}
\mL(\cG) = 
\begin{bmatrix*}[r]
    &2&\text{-}1&\text{-}1 \\ 
    &\text{-}1&1&0 \\ 
    &\text{-}1&0&1 \\ 
\end{bmatrix*}.
\end{eqnarray*}
\pars{Note, that the} \ar{condition number of $\mW$ is better} \pars{as compared to the one of} \ar{$\mL(\cG)$: we have $\chi(\mW) = 1$ and $\chi(\mL(\cG)) = 3$.}

\subsection{Spectrum of $\mW$}\label{subsec:complex_lapl_spectrum}

In this section, we analyze the spectrum of the novel communication matrix $\mW$.
\begin{lemma}\label{lemma:complex_lapl_spectrum}
    For matrix $\mW$ defined in \eqref{eq:complex_lapl} we have $\bLambda(\mW) = \bigcup\limits_{\ell=1}^n \bLambda(\mL(\cG^\ell))$.
\end{lemma}
\begin{pf}
Firstly, let $\mW\bx = \lambda\bx$ for some $\bx\ne 0$ and $\lambda\in\mathbb{C}$. By definition of $\mW$, we have $\mL(\cG^\ell) \bx^\ell = \lambda \bx^\ell$ for $\ell = 1, \ldots, n$. Since $\bx\ne 0$, there exists $\ell$ such that $\bx^\ell\neq 0$, and therefore $\lambda\in \bigcup\limits_{\ell=1}^n \bLambda(\mL(\cG^\ell))$.

Secondly, let $\mL(\cG^\ell) \bx^\ell = \lambda \bx^\ell$ for some $\ell = 1, \ldots, n$. Setting $\bx = \bx^\ell\otimes {\be_n^{(\ell)}}$ we obtain 
\begin{align*}
\mW\bx &= \sum_{j=1}^n (\mL(\cG^\ell)\otimes \be_n^{(j)}\be_n^{(j)\top})(\bx^\ell \otimes {\be_n^{(\ell)}}) \\
&= (\mL(\cG^\ell) \bx^\ell) \otimes {\be_n^{(\ell)}} = \lambda(\bx^\ell\otimes{\be_n^{(\ell)}}) = \lambda\bx
\end{align*}
i.e. $\lambda\in\bLambda(\mW)$. As a result, $\bLambda(\mW)\subseteq\bigcup\limits_{\ell=1}^n \bLambda(\mL(\cG^\ell))$ and $\bigcup\limits_{\ell=1}^n \bLambda(\mL(\cG^\ell))\subseteq \bLambda(\mW)$, therefore, $\bigcup\limits_{\ell=1}^n \bLambda(\mL(\cG^\ell)) = \bLambda(\mW)$.
\end{pf}

It immediately follows from Lemma \ref{lemma:complex_lapl_spectrum} that $\ds\lambda_{\max}(\mW) = \max_{1\leq\ell\leq n} \lambda_{\max}(\mL(\cG^\ell))$ and $\ds\lambda_{\min}^+(\mW) = \min_{1\leq\ell\leq n} \lambda_{\min}^+(\mL(\cG^\ell))$. Thus, we have
\begin{equation}\label{chi}
\chi(\mW) = \frac{\underset{1\leq\ell\leq n}{\max} \lambda_{\max}(\mL(\cG^\ell))}{\underset{1\leq\ell\leq n}{\min} \lambda_{\min}^+(\mL(\cG^\ell))}.
\end{equation}
\ar{Concerning the example in Figure \ref{fig:three graphs}, we have $\lambda_{\max}(\blue{\mL(\cG^1)}) = \lambda_{\min}^+(\blue{\mL(\cG^1)}) = 1$, $\lambda_{\max}(\green{\mL(\cG^2)}) = \lambda_{\min}^+(\green{\mL(\cG^2)}) = 1$ and therefore $\chi(\mW) = 1$.
}
\begin{remark}
We can determine $\mL(\cG)$ up to a positive multiplicative constant. Since that we can consider $\lambda_{\max}(\mL(\cG^\ell))= 1$ for all $l=1,...,n$ without loss of generality. For that we need some preprocessing to estimate $\left\{\lambda_{\max}(\mL(\cG^\ell))\right\}_{l=1}^n$ by using Power method (see \cite{golub2013matrix}) that could be done in a distributed manner. We also note that $\chi(\mW)$ from \eqref{chi} can be bounded from below by the largest diameter of graphs $\cG^\ell$, $l=1,...,n$. Moreover, for many important classes of graphs this lower bound is tight up to a $\ln n$ factor \cite{scaman2017optimal}.
\end{remark}




%% file: sections/ring_of_cliques.tex
\section{Illustrative example: cycle of cliques}\label{sec:ring_of_cliques}

We illustrate the effect of using matrix $\mW$ defined in \eqref{eq:complex_lapl} on a synthetic example. \pars{Our test case is a computing network with a two-level hierarchical structure. In the lower layer we have a clique (i.e. a fully connected graph) with $k\geq2$ nodes. Each node in a clique is supposed to share a local variable with its neighbors. In the upper level $n$ cliques communicate through undirected cyclic topology: every clique has a ``negotiator'', i.e. a node that has links with similar neighboring nodes (first node of every clique w.l.o.g.).} We denote our graph $\cG_{RC}$.

\begin{figure}[h!]
    \centering
    \hspace{0.5cm}
    \begin{subfigure}[t]{0.3\linewidth}
        \centering
        \includegraphics[width=1.0\linewidth]{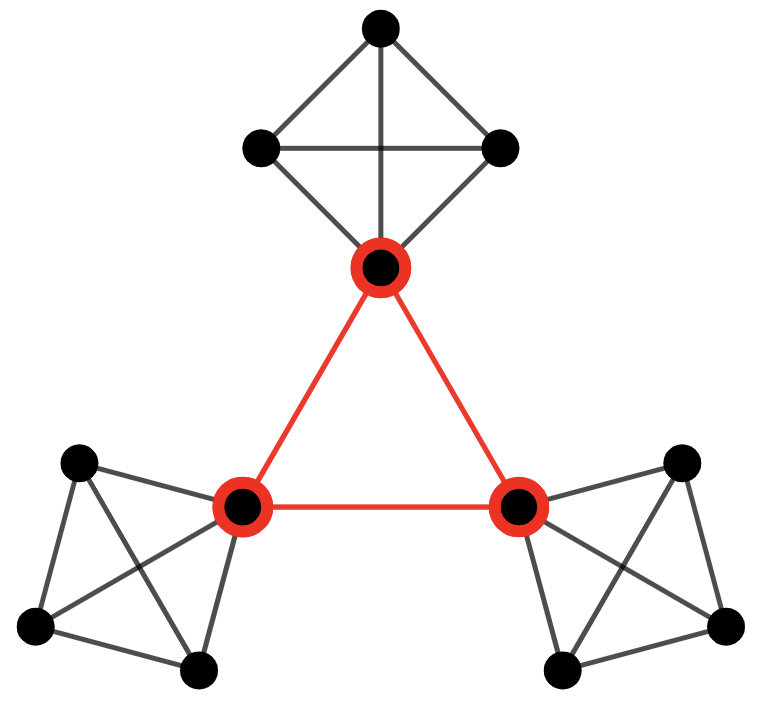}
        \caption{Graph $\cG_{RC}$ with $n=3$ cliques of size $k=4$}
        \label{fig:ring_clique_example}
    \end{subfigure}
    \hfill
    \hspace{0.5cm}
    \begin{subfigure}[t]{0.25\textwidth}
        \centering
        \includegraphics[width=0.35\linewidth]{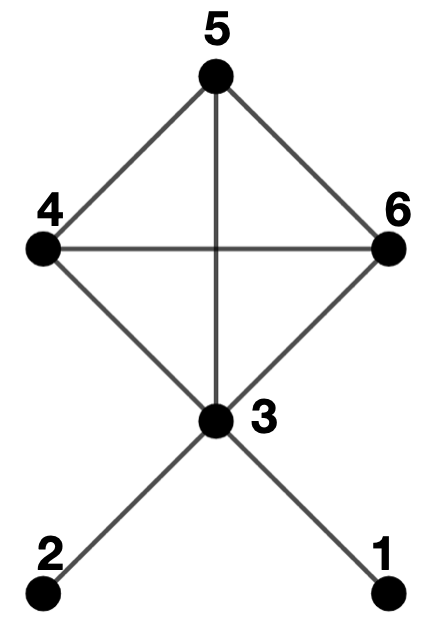}
        \caption{Subgraph $\cG_{cr}$ corresponding to \\variable $x^1$}
        \label{fig:4crown}
    \end{subfigure}
    \caption{Hierarchical graph with cliques and its "crown" subgraph}
    \centering
\end{figure}

\pd{For every clique, we consider its union with two adjacent nodes and call such subgraph a "crown", see Figure \ref{fig:4crown}. Overall, we have $n$ "crowns" $\tilde\cG^\ell$, $\ell = 1,\ldots,n$, each corresponding to one of the $n$ variables, i.e., each function $\tilde{f}_i$ held by a vertex $i$ of $\tilde\cG^\ell$ depends on $x^\ell$. Further, if a node $i$ lies in the intersection of two "crowns", its part $\tilde{f}_i$ of the objective depends also on the variables corresponding to each of the neighboring "crowns".} \ar{All "crowns" $\tilde\cG^\ell$ are isomorphic and further denote the "crown" graph as $\cG_{cr}$.}


Consider an example of graph $\cG_{RC}$ with $n=3$ given in Figure \ref{fig:ring_clique_example} and enumerate the vertices of the top "crown" as shown in Figure \ref{fig:4crown}. \pars{Here, the lower level of hierarchy are the fully connected cliques, depicted in black. The upper level of communication is an undirected ring (in red). The ``crowns'' can be treated as the cliques supplemented by two extra vertices and corresponding edges}. Then the functions $\tilde{f}_4, \tilde{f}_5, \tilde{f}_6$ depend only on the variable $x^1$  since the corresponding nodes belong only to the top "crown". On the other hand, each of the nodes 1, 2, 3 lies in the intersection of all 3 crowns and therefore functions $\tilde{f}_1, \tilde{f}_2, \tilde{f}_3$ depend on $x^1, x^2, x^3$.
In this section, we 
 illustrate that matrix $\mW$ defined in \eqref{eq:complex_lapl} has a better condition number than the Laplacian $\mL(\cG_{RC})$. 
\begin{theorem}\label{th:ring_clique_condition_number}
For the hierarchical ring-clique graph it holds $$
\chi(\mW) = \Theta(k),~ \chi(\mL(\cG_{RC})) = \Theta(n^2k^2).
$$
\end{theorem}
\ar{We prove Theorem \ref{th:ring_clique_condition_number} in a sequence of Lemmas \ref{lem:crown-graph_cond_number}, \ref{lem:lapl_spectrum_union}, \ref{lem:matrix_determ}, \ref{lem:estimate_lambda_min} presented below in this section.} The full proofs of the Lemmas are presented in Appendix.
\pd{\ar{Theorem \ref{th:ring_clique_condition_number}} illustrates the flexibility and efficiency of our approach compared to standard approaches that do not take into account the partitioned structure of problem \eqref{eq:problem_initial}.} \ar{Substituting matrix $\mW$ instead of $\mL(\cG_{RC})$ allows to enhance the convergence rate of distributed algorithms. In the following corollary, we illustrate this speedup on state-of-the-art primal and dual optimization methods.}
\ar{
\begin{corollary}
    Let all the functions $f_i$, be $L$-smooth, $\mu$ strongly-convex and stored at the nodes of $\cG_{RC}$. Consider two algorithms: dual MSDA \cite{scaman2017optimal} and primal OPAPC \cite{kovalev2020optimal}. If we use $\mL(\cG_{RC})$, each of these methods has communication complexity $O\cbraces{\sqrt{\frac{L}{\mu}}nk \ln\frac{1}{\eps}}$; the communication complexity becomes $O\cbraces{\sqrt{\frac{L}{\mu}} \sqrt k\ln\frac{1}{\eps}}$ if we use $\mW$, where $\eps$ is the accuracy.
\end{corollary}
}
\pars{\begin{remark}
Obviously, in general case there may be examples of graphs where our approach does not provide significant improvement in terms of condition numbers of the corresponding Laplacian matrices. However, the essential feature of partitioned representation of the state vector with sharing of the necessary states only makes this formulation of the distributed optimization problem attractive in terms of more sparse communication topology and reduced information exchange (as compared to \eqref{eq:problem_affine_contraints}). It also corresponds to the preservation of the privacy of the data of interacting groups of nodes. 
\end{remark}}

We first estimate $\chi(\mW)$. From Lemma \ref{lemma:complex_lapl_spectrum} it follows that $\bLambda(\mW) = \bLambda(\mL(\cG_{cr}))$, where $\cG_{cr}$ denotes the "crown" graph. To estimate the asymptotics of the condition number  of "crown" graph, i.e., the condition number of its Laplacian,  we use a technique described in \cite{pozrikidis2014introduction}.
\begin{lemma}\label{lem:crown-graph_cond_number}
The condition number of crown-graph has asymptotics $\chi(\mL(\cG_{cr})) = \Theta(k)$, where $k$ is the clique size.
\end{lemma}
\begin{wrapfigure}[12]{r}{2.5cm}
\vspace{-0.3cm}
\includegraphics[width=2.5cm]{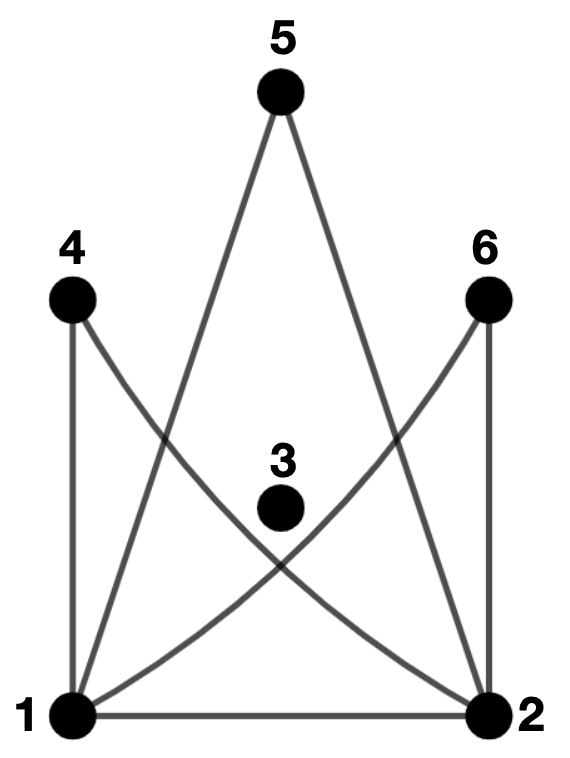}
\caption{Complement of $\cG_{rc}$}
\label{fig:crown_complement}
\end{wrapfigure}
\begin{pf}
Firstly, to find the eigenvalues of \pd{$\mL(\cG_{cr})$, we find the eigenvalues of $\mL(\bar{\cG}_{cr})$, where $\bar{\cG}_{cr}$ is the  complement of $\cG_{cr}$.}

    The complement has one isolated vertex (node 3 in Figure \ref{fig:crown_complement}). We denote the Laplacian of the connected part of $\bar\cG_{cr}$ (i.e. graph in Figure \ref{fig:crown_complement} with vertices $\braces{1, 2, 4, 5, 6}$) as $\mL'(\bar{\cG}_{cr})$.
    \vspace{0.3cm}
    \begin{align*}
        \mL'&(\bar\cG_{cr}) = \diag(k + 2, k + 2, 2, \ldots, 2) \\
        &- \one_{k+1}(\be_{k+1}^{(1)} + \be_{k+1}^{(2)})^\top - (\be_{k+1}^{(1)} + \be_{k+1}^{(2)}) \one_{k+1}^\top
    \end{align*}
    
    Eigenvalues $\lambda'$ of $\mL'(\bar\cG_{cr})$ are defined through the equation
    \begin{align*}
        \text{det} \left( \mL' \left( \bar\cG_{cr} \right) - \lambda' \mI_{k+1} \right) = 0
    \end{align*}
    which, via linear conversions (see Appendix~\ref{app:crown-graph_cond_number} for details), leads us to the equation 
    \begin{align*}
    \lambda' (\lambda'-k-1)^2(\lambda'-2)^{k-2} = 0.
    \end{align*}
    Thus, the eigenvalues of \pd{$\mL'\left(\bar{\cG}_{cr}\right)$} have the form $\lambda_1' = 0,~ \lambda_2' = \lambda_3' = k + 1,~ \lambda_4' = \ldots = \lambda_{k+1}' = 2$. Using Eq. $2.2.22$ from \cite{pozrikidis2014introduction}, we obtain the eigenvalues of \pd{$\mL\left(\cG_{cr}\right)$} to be $\lambda_1 = 0,~ \lambda_2 = k + 2,~ \lambda_3 = \lambda_4 = 1,~ \lambda_5 = \ldots = \lambda_{k+2} = k$. Thus, we have $\chi(\mL(\cG_{cr})) = (k+2)/1 = \Theta(k)$.

        
        
        
        
\end{pf}

Applying Lemma~\ref{lemma:complex_lapl_spectrum} we obtain $\chi(\mW) = \chi(\mL(\cG_{cr})) = \Theta(k)$.

Our next goal is to estimate $\chi(\mL(\pd{\cG_{RC}}))$ and show that it is worse than $\chi(\mW)$.
To compute $\chi(\mL(\cG_{RC}))$, we decompose $\mL(\cG_{RC})$ into Laplacians of a $k$-clique $\cG_{C}$ and a ring graph $\cG_{R}$ that has $n$ nodes. We also introduce matrix $\mB = \be_k^{(1)}\be_k^{(1)\top}$
that allows us to write $\mL(\cG_{RC})$ in the following form.
\begin{eqnarray}\label{eq:lapl_ring_of_cliques}
    \bf{\mL(\cG_{RC})} = \mI_n \otimes \mL(\cG_{C}) + \mL(\cG_{R}) \otimes \mB
\end{eqnarray}
\ar{Let $\braces{\lambda_i}_{i=1}^n$ be the eigenvalues of $\cG_R$}. To obtain the eigenvalues of $\mL(\cG_{RC})$, we diagonalize it and decompose its spectrum. The result is formulated in the following Lemma.
\begin{lemma}\label{lem:lapl_spectrum_union}
$\bLambda(\mL(\cG_{RC})) = \bigcup\limits_{i=1}^n \bLambda(\mL(\cG_{C}) + \lambda_i\mB)$. 
\end{lemma}
\begin{pf}
Let $\mL(\cG_{R}) = \mS \Psi \mS^\top$, where $\Psi = \diag\left(\lambda_1, \lambda_2, \ldots, \lambda_n \right)$ is a diagonal matrix with eigenvalues of $\mL(\cG_{R})$ and $\mS\mS^\top = \mI_n$. We define
\begin{eqnarray*}
    \hat\mL = \left( \mS \otimes \mI_k \right)^{\top} \mL(\cG_{RC}) \left( \mS \otimes \mI_k \right) = \mI_n \otimes \mL(\cG_{C}) + \Psi \otimes \mB.
\end{eqnarray*}
We have $\bLambda(\hat\mL) = \bLambda(\mL(\cG_{RC}))$. Indeed, let $\bx$ be an eigenvector of $\mL(\cG_{RC})$ such that $\mL(\cG_{RC})\bx = \theta\bx$. Then
\begin{align*}
    \hat\mL\cdot ((\mS\otimes\mI_k)^{\top}\bx) = \theta\cdot ((\mS\otimes\mI_k)^{\top}\bx),
\end{align*}
i.e., $(\mS\otimes\mI_k)^\top\bx$ is an eigenvector of $\hat\mL$.

Further, for any eigenvalue $\theta$ of $\hat\mL$, we have
\begin{eqnarray*}
\hat\mL 
\begin{bmatrix}
x_1 \\
\vdots \\
x_n
\end{bmatrix}
=
\begin{bmatrix}
 \left(\mL(\cG_{C}) + \lambda_1 \mB \right) x_1  \\
\vdots  \\
\left(\mL(\cG_{C}) + \lambda_n \mB \right) x_n  
\end{bmatrix}
= \theta
\begin{bmatrix}
    x_1 \\
    \vdots \\
    x_n
\end{bmatrix},
\end{eqnarray*}
i.e., $(\mL(\cG_{C}) + \lambda_i\mB)x_i = \theta x_i$ for $i = 1, \ldots, n$. Consequently, $\bLambda(\hat\mL) = \bigcup\limits_{i=1}^n \bLambda(\mL(\cG_{C}) + \lambda_i\mB)$, which concludes the proof.
\end{pf}
Due to Lemma \ref{lem:lapl_spectrum_union}, we only have to find the spectra $\bLambda(\mL(\cG_C) + \lambda_i\mB)$ for $i = 1, \ldots, n$. This is done using the matrix determinant lemma.
\begin{lemma}\label{lem:matrix_determ}
The matrix $\mL(\cG_C) + \lambda\mB$ has the following eigenvalues.
\begin{subequations}\label{eq:subspectrum}
\begin{align}
    \theta_1(\lambda) &= \dfrac{k+\lambda + \sqrt{\lambda^2 + 2(k-2)\lambda+k^2}}{2}, \label{eq:subspectrum_1} \\
    \theta_2(\lambda) &= \dfrac{k+\lambda - \sqrt{\lambda^2 + 2(k-2)\lambda+k^2}}{2}, \label{eq:subspectrum_2}\\
    \theta_3(\lambda) &= \ldots = \theta_n(\lambda) = k. \nonumber
\end{align}
\end{subequations}
\end{lemma}
\begin{pf}
Let $\theta$ denote some eigenvalue of $\mL(\cG_C) + \lambda\mB$. Then, $\mL(\cG_{C}) + \lambda\mB - \theta\mI_k = \mA - \one_k\one_k^\top$, where $\mA = \diag \left(k+\lambda-\theta, k-\theta, \ldots, k-\theta \right)$.

Firstly, note that $\theta = k$ is an eigenvalue of $\mL(\cG_C) + \lambda\mB$. Secondly, $\theta = k + \lambda$ is an eigenvalue of $\mL(\cG_C) + \lambda\mB$ if and only if $\lambda = 0$ (details can be found in Appendix~\ref{app:matrix_determ}).

Further, we assume that $\theta\ne k$ and $\theta\ne k + \lambda$. According to the matrix determinant lemma we have 
\begin{align*}
    \text{det}& \left( \mA - \one_k\one_k^\top \right) = \left( 1 - \one^\top \cdot \mA^{-1} \cdot \one \right) \cdot \text{det} \mA \\
    &= (k - \theta)^{k-2} (\theta - \theta_1) (\theta - \theta_2)
\end{align*}
where $\theta_1$ and $\theta_2$ are defined in \eqref{eq:subspectrum}. Details can be found in Appendix~\ref{app:matrix_determ}.
\end{pf}
Equation \eqref{eq:subspectrum} gives explicit formulas for eigenvalues of $\mL(\cG_{RC})$, and it remains to estimate $\lambda_{\max}(\mL(\cG_{RC}))$ and $\lambda_{\min}^+(\mL(\cG_{RC}))$. The eigenvalues of the ring graph $\cG_R$ have the form $\lambda_i = 2 - 2\cos\frac{2\pi i}{n}$. Note that $0\leq \lambda_i\leq 4$.


For the largest eigenvalue of $\mL(\cG_{RC})$, we have
$$
\lambda_{\max}(\mL(\cG_{RC})) \leq \frac{k + 4 + \sqrt{4 + 2(k - 2)\cdot 2 + k^2}}{2} = \Theta(k).
$$
Estimating $\lambda_{\min}^+(\mL(\cG_{RC}))$ is less straightforward and we need to compute the minimal $\theta_2$ defined in \eqref{eq:subspectrum_2} over $\lambda\in\braces{\lambda_1, \ldots, \lambda_n}$.
\begin{lemma}\label{lem:estimate_lambda_min}
It holds that $\lambda_{\min}^+(\mL(\cG_{RC})) = \Theta(\frac{1}{n^2 k})$.
\end{lemma}
\begin{pf}
Firstly, let us show that $\theta_2(\lambda)$ defined in \eqref{eq:subspectrum_2} is monotonically increasing by considering its derivative:
\begin{eqnarray*}
    \frac{d\theta_2(\lambda)}{d\lambda} = \dfrac{-k-\lambda+2}{2\sqrt{k^2 + \lambda^2 + 2(k-2)\lambda}} + \dfrac{1}{2}
\end{eqnarray*}
Since $k \geq 2$ it holds $k+\lambda - 2 < \sqrt{k^2 + \lambda^2 + 2(k-2)\lambda}$ and we have $\frac{d\theta_2(\lambda)}{d\lambda} > 0$. Since $\theta_2(0) = 0$, the minimal positive eigenvalue of $\mL(\cG_{RC})$ is reached at minimal positive $\lambda$, that is, $\lambda_{\min}^+(\cG_{RC}) = \theta_2(2 - 2\cos\frac{2\pi}{n})$.

Secondly, we approximate $\theta_2(2 - 2\cos\frac{2\pi}{n})$ using the Taylor series at $n\to\infty$ (see Appendix~\ref{app:estimate_lambda_min} for details) and get
\begin{align*}
    \theta_2\cbraces{2 - 2\cos\frac{2\pi}{n}} = \frac{4\pi^2}{n^2k} + o\cbraces{\frac{1}{n^3}}.
\end{align*}
It follows that $\lambda_{\min}^+(\cG_{RC}) = \Theta(\frac{1}{n^2k})$. 
\end{pf}










    
    

%% file: sections/conclusion.tex
\section{Conclusion}

In this paper, we described a simple and universal framework that allows to work with distributed optimization problems with partitioned structure. The proposed method enables to write communication constraints as affine constraints by using a modified Lagrangian matrix. This reformulation lowers the cost of each communication round and also reduces the communication complexities of distributed optimization algorithms.

%% file: sections/appendix.tex
{\begin{center}\large \textbf{Appendix} \end{center}}

Let $\mA$ and $\mB$ be quadratic matrices of the same size. Throughout the Appendix, we equivalently write $|\mA|$ and $\det(\mA)$ for determinant of $\mA$. We also write $\mA\sim\mB$ if $\mA$ is obtained from $\mB$ by adding one column (row) multiplied by a scalar to another column (row). Note that if $\mA\sim\mB$ then $\det(\mA) = \det(\mB)$.

\section{Proof of Lemma \ref{lem:crown-graph_cond_number}}\label{app:crown-graph_cond_number}

Firstly, to find the eigenvalues of \pd{$\mL(\cG_{cr})$, we find the eigenvalues of $\mL(\bar{\cG}_{cr})$, where $\bar{\cG}_{cr}$ is the  complement of $\cG_{cr}$.}

The complement has one isolated vertex (node 3 in Figure \ref{fig:crown_complement}). We denote the Laplacian of the connected part of $\bar\cG_{cr}$ (i.e. graph in Figure \ref{fig:crown_complement} with vertices $\braces{1, 2, 4, 5, 6}$) as $\mL'(\bar{\cG}_{cr})$.
\begin{align*}
\mL'&(\bar\cG_{cr}) = \diag(k + 2, k + 2, 2, \ldots, 2) - \one_{k+1}(\be_{k+1}^{(1)} + \be_{k+1}^{(2)})^\top - (\be_{k+1}^{(1)} + \be_{k+1}^{(2)}) \one_{k+1}^\top.
\end{align*}

Eigenvalues $\lambda'$ of $\mL'(\bar\cG_{cr})$ are defined through the equation
\begin{align*}
\text{det} \left( \mL' \left( \bar\cG_{cr} \right) - \lambda' \mI_{k+1} \right) = 0.
\end{align*}
Now we find the roots of the equation above.

\begin{align*}
    \mL'(\ol\cG_{cr}) - \lambda'\mI_{k+1}
    &= 
    \begin{bmatrix}
    &k-\lambda'~~&-1~~&-1~~&\cdots&\cdots&-1 \\
    &-1~~&k-\lambda'~~&-1~~&\cd&\cd&-1 \\
    &-1~~&-1~~&2-\lambda'~~&0&\cd&0 \\
    &\vd~~&\vd~~&0~~&\ddots&&\vd \\
    &\vd~~&\vd~~&\vd~~&&\ddots&0 \\
    &-1~~&-1~~&0~~&\cd&0&2-\lambda' \\
    \end{bmatrix}
    \sim
    \begin{bmatrix}
    &k-\lambda'+1~~&-1~~&-1~~&\cdots&\cdots&-1 \\
    &\lambda'-k-1~~&k-\lambda'~~&-1~~&\cd&\cd&-1 \\
    &0~~&-1~~&2-\lambda'~~&0&\cd&0 \\
    &\vd~~&\vd~~&0~~&\ddots&&\vd \\
    &\vd~~&\vd~~&\vd~~&&\ddots&0 \\
    &0~~&-1~~&0~~&\cd&0&2-\lambda' \\
    \end{bmatrix} \\
    &\sim
    \begin{bmatrix}
    &2k-2\lambda'+2~~&\lambda'-k-1~~&0~~&\cdots&\cdots&0 \\
    &\lambda'-k-1~~&k-\lambda'~~&-1~~&\cd&\cd&-1 \\
    &0~~&-1~~&2-\lambda'~~&0&\cd&0 \\
    &\vd~~&\vd~~&0~~&\ddots&&\vd \\
    &\vd~~&\vd~~&\vd~~&&\ddots&0 \\
    &0~~&-1~~&0~~&\cd&0&2-\lambda' \\
    \end{bmatrix}.
\end{align*}

Then
\begin{align*}
&\text{det}
\begin{bmatrix}
    &k-\lambda'~~&-1~~&-1~~&\cdots&\cdots&-1 \\
    &-1~~&k-\lambda'~~&-1~~&\cd&\cd&-1 \\
    &-1~~&-1~~&2-\lambda'~~&0&\cd&0 \\
    &\vd~~&\vd~~&0~~&\ddots&&\vd \\
    &\vd~~&\vd~~&\vd~~&&\ddots&0 \\
    &-1~~&-1~~&0~~&\cd&0&2-\lambda' \\
    \end{bmatrix} \\
    &\qquad=
    (2k-2\lambda'+2)
    \begin{vmatrix}
    &k-\lambda'~~&-1~~&\cd&\cd&-1 \\
    &-1~~&2-\lambda'~~&0&\cd&0 \\
    &\vd~~&0~~&\ddots&&\vd \\
    &\vd~~&\vd~~&&\ddots&0 \\
    &-1~~&0~~&\cd&0&2-\lambda' \\   
    \end{vmatrix}
     - (\lambda'-k-1)
     \begin{vmatrix}
    &\lambda'-k-1~~&-1~~&\cd&\cd&-1 \\
    &0~~&2-\lambda'~~&0&\cd&0 \\
    &\vd~~&0~~&\ddots&&\vd \\
    &\vd~~&\vd~~&&\ddots&0 \\
    &0~~&0~~&\cd&0&2-\lambda' \\     
    \end{vmatrix}.
\end{align*}

For the second term we have
\begin{eqnarray*}
    \begin{vmatrix}
    &\lambda'-k-1~~&-1~~&\cd&\cd&-1 \\
    &0~~&2-\lambda'~~&0&\cd&0 \\
    &\vd~~&0~~&\ddots&&\vd \\
    &\vd~~&\vd~~&&\ddots&0 \\
    &0~~&0~~&\cd&0&2-\lambda' \\     
     \end{vmatrix} = (\lambda'-k-1)(2-\lambda')^{k-1}.
\end{eqnarray*}

For the first term we have
\begin{align*}
&
\begin{vmatrix}
    &k-\lambda'~~&-1~~&\cd&\cd&-1 \\
    &-1~~&2-\lambda'~~&0&\cd&0 \\
    &\vd~~&0~~&\ddots&&\vd \\
    &\vd~~&\vd~~&&\ddots&0 \\
    &-1~~&0~~&\cd&0&2-\lambda' \\   
    \end{vmatrix} \\
	&\qquad=
    (k-\lambda')(2-\lambda')^{k\text{-}1}
    -
    (-1) \cdot 
    \begin{vmatrix}
     &-1~~&0~~&\cdots&0 \\
     &-1~~&2-\lambda'~~&\cdots&0 \\
     &\vdots~~&\vdots~~&\ddots&\vdots \\
     &-1~~&0~~&\cdots&2-\lambda' \\
    \end{vmatrix}
    +
    (-1) \cdot
    \begin{vmatrix}
     &-1~~&2-\lambda'~~&0&\cdots&0 \\
     &-1~~&0~~&0&\cdots&0 \\
     &-1~~&0~~&2-\lambda'&&\vdots \\
     &\vdots~~&\vdots~~&&\ddots&0 \\
     &-1~~&0~~&\cdots&0&2-\lambda' \\
    \end{vmatrix}
    -(-1) \cdot \ldots \\\nonumber
    &\qquad=
    (k-\lambda') (2-\lambda')^{k-1} - (k-1)\cdot(2-\lambda')^{k-2}.
\end{align*}

Consequently,
\begin{align*}
    \det(\mL'(\ol\cG_{cr}) - \lambda'\mI_{k+1}) &= (2k-2\lambda'+2)\left( (k-\lambda')(2-\lambda')^{k-1}-(k-1)(2-\lambda')^{k-2} \right) - (\lambda' - k - 1)^2 (2-\lambda')^{k-1} \\
    &= \left[ (k-\lambda'+1)(2-\lambda')^{k-2} \right] \cdot \left[ 2(k-\lambda')(2-\lambda') - 2(k-1) - (k-\lambda'-1)(2-\lambda') \right] \\
    &= (k-\lambda'+1)(2-\lambda')^{k-2}\cdot [\lambda'(\lambda'-k-1)].
\end{align*}
which leads us to the equation 
\begin{align*}
    \lambda' (\lambda'-k-1)^2(\lambda'-2)^{k-2} = 0.
\end{align*}

Thus, the eigenvalues of \pd{$\mL'\left(\bar{\cG}_{cr}\right)$} have the form $\lambda_1' = 0,~ \lambda_2' = \lambda_3' = k + 1,~ \lambda_4' = \ldots = \lambda_{k+1}' = 2$. Using Eq. $2.2.22$ from \cite{pozrikidis2014introduction}, we obtain the eigenvalues of \pd{$\mL\left(\cG_{cr}\right)$} to be $\lambda_1 = 0,~ \lambda_2 = k + 2,~ \lambda_3 = \lambda_4 = 1,~ \lambda_5 = \ldots = \lambda_{k+2} = k$. Thus, we have $\chi(\mL(\cG_{cr})) = (k+2)/1 = \Theta(k)$.

\section{Proof of Lemma~\ref{lem:matrix_determ}}\label{app:matrix_determ}

Let $\theta$ denote some eigenvalue of $\mL(\cG_C) + \lambda\mB$. Then, $\mL(\cG_{C}) + \lambda\mB - \theta\mI_k = \mA - \one_k\one_k^\top$, where
\begin{eqnarray*}
	\mA = \diag \left(k+\lambda-\theta, k-\theta, \ldots, k-\theta \right) = 
	\begin{bmatrix}
		&k+\lambda-\theta&0&\cd&\cd&0 \\
		&0&k-\theta&0&\cd&0 \\
		&\vd&0&\ddots&\ddots&\vd \\
		&\vd&\vd&\ddots&\ddots&0 \\
		&0&0&\cd&0&k-\theta \\
	\end{bmatrix}
\end{eqnarray*}

Firstly, note that $\theta = k$ is an eigenvalue of $\mL(\cG_C) + \lambda\mB$. 

Secondly, note that $\theta = k + \lambda$ is an eigenvalue of $\mL(\cG_C) + \lambda\mB$ if and only if $\lambda = 0$. Indeed, for $\theta = k + \lambda$ it holds
\begin{align*}
	\mA - \one_k\one_k^\top &= 
	\begin{bmatrix}
		-1~~& -1~~& -1~~& \ldots~~& -1 \\
		-1~~& -\lambda-1~~& -1~~& \ldots~~& -1 \\
		-1~~& -1~~& -\lambda-1~~& \ldots~~& -1 \\
		\vdots~~& \vdots~~& \vdots& \ddots~~& ~~\vdots \\
		-1~~& -1~~& -1~~& \ldots~~& -\lambda - 1
	\end{bmatrix}
	\sim
	\begin{bmatrix}
	-1~~& -1~~& -1~~& \ldots~~& -1 \\
	~~0~~& -\lambda~~& ~~0~~& \ldots~~& ~~0 \\
	-1~~& -1~~& -\lambda-1~~& \ldots~~& -1 \\
	\vdots~~& \vdots~~& \vdots& \ddots~~& ~~\vdots \\
	-1~~& -1~~& -1~~& \ldots~~& -\lambda - 1
	\end{bmatrix} \\
	&\sim
	\begin{bmatrix}
	-1~~& 0~~& -1~~& \ldots~~& -1 \\
	~~0~~& -\lambda~~& ~~0~~& \ldots~~& ~~0 \\
	-1~~& 0~~& -\lambda-1~~& \ldots~~& -1 \\
	\vdots~~& \vdots~~& \vdots& \ddots~~& ~~\vdots \\
	-1~~& 0~~& -1~~& \ldots~~& -\lambda - 1
	\end{bmatrix}
	\sim\ldots\sim
	\begin{bmatrix}
	-1~~& ~~0~~& ~~0~~& \ldots~~& ~~0 \\
	~~0~~& -\lambda~~& ~~0~~& \ldots~~& ~~0 \\
	~~0~~& ~~0~~& -\lambda~~& \ldots~~& ~~0 \\
	\vdots~~& \vdots~~& \vdots& \ddots~~& ~~\vdots \\
	~~0~~& ~~0~~& ~~0~~& \ldots~~& -\lambda
	\end{bmatrix}.
\end{align*}
Therefore, $\det(\mA - \one_k\one_k^\top) = (-1)^k\lambda^{k-1}$ and $\det(\mA - \one_k\one_k^\top)$ if and only if $\lambda = 0$. We conclude that either $\theta = k + \lambda = k$ (when $\lambda = 0$) or $\theta = k + \lambda$ is not an eigenvalue of $\mL(\cG_C + \lambda\mB)$ (when $\lambda\neq 0$).

Further, we assume that $\theta\ne k$ and $\theta\ne k + \lambda$.

By definition of $\mA$ we have $\det(\mA) = (k + \lambda - \theta) (k - \theta)^{k-1}\ne 0$.

According to matrix determinant lemma (see Lemma $1.1$ from \cite{ding2007eigenvalues}):
\begin{eqnarray*}
	\text{det} \left( \mA - \one_k \one_k^\top \right) = \left( 1 - \one_k^\top \mA^{-1}\one_k \right) \cdot \text{det} \mA
\end{eqnarray*}

Since $\theta\neq k$ and $\theta\neq k + \lambda$ we have
\begin{eqnarray*}
	\mA^{-1} = 
	\begin{bmatrix}
		&\dfrac{1}{k+\lambda-\theta}&0&\cd&\cd&0 \\
		&0&\dfrac{1}{k-\theta}&0&\cd&0 \\
		&\vd&0&\ddots&\ddots&\vd \\
		&\vd&\vd&\ddots&\ddots&0 \\
		&0&0&\cd&0&\dfrac{1}{k-\theta} \\
	\end{bmatrix}
\end{eqnarray*}

It follows that $1 - \one_k^\top\mA^{-1}\one_k = 1 -\dfrac{1}{k+\lambda-\theta} - \dfrac{k-1}{k-\theta}$.

\begin{eqnarray*}
	\text{det}(\mL(\cG_C) + \lambda\mB - \theta\mI_k) &=& \left( 1 -\dfrac{1}{k+\lambda-\theta} - \dfrac{k-1}{k-\theta} \right) (k + \lambda - \theta) (k-\theta)^{k-1}  \nonumber \\
	&&= (k+\lambda-\theta)(k-\theta)^{k-1} - (k-\theta)^{k-1} - (k-1)(k+\lambda-\theta) (k-\theta)^{k-2} \nonumber \\
	&&= (k-\theta)^{k-2} \cdot \left( (k + \lambda - \theta) \cdot (k-\theta) - (k-\theta) - (k-1)(k+\lambda-\theta) \right) \nonumber \\
	&&= (k-\theta)^{k-2} \left( \theta^2 + \theta(-k-\lambda) + \lambda \right) \nonumber \\
	&&= (k-\theta)^{k-2} \left( \theta - \theta_1 \right)\left( \theta - \theta_2 \right).
\end{eqnarray*}
Here
\begin{align*}
\theta_1 &= \dfrac{k+\lambda + \sqrt{\lambda^2 + 2(k-2)\lambda+k^2}}{2}, \\
\theta_2 &= \dfrac{k+\lambda - \sqrt{\lambda^2 + 2(k-2)\lambda+k^2}}{2}.
\end{align*}
The rest eigenvalues are $\theta_3 = \ldots = \theta_k = k$.

\section{Proof of Lemma~\ref{lem:estimate_lambda_min}}\label{app:estimate_lambda_min}

Firstly, let us show that $\theta_2(\lambda)$ defined in \eqref{eq:subspectrum_2} is monotonically increasing by considering its derivative:
\begin{eqnarray*}
	\frac{d\theta_2(\lambda)}{d\lambda} = \dfrac{-k-\lambda+2}{2\sqrt{k^2 + \lambda^2 + 2(k-2)\lambda}} + \dfrac{1}{2}
\end{eqnarray*}
Since $k \geq 2$, it holds $k+\lambda - 2 < \sqrt{k^2 + \lambda^2 + 2(k-2)\lambda}$ and we have $\frac{d\theta_2(\lambda)}{d\lambda} > 0$. Since $\theta_2(0) = 0$, the minimal positive eigenvalue of $\mL(\cG_{RC})$ is reached at minimal positive $\lambda$, that is, $\lambda_{\min}^+(\cG_{RC}) = \theta_2(2 - 2\cos\frac{2\pi}{n})$.

Secondly, we approximate $\theta_2(2 - 2\cos\frac{2\pi}{n})$ using the Taylor series at $n\to\infty$.
\begin{align*}
	\theta_2 \cbraces{2 - 2\cos\frac{2\pi}{n}}
	&= \frac{1}{2}\sbraces{k + 2\cbraces{\frac{4\pi^2}{2n^2} + o\cbraces{\frac{1}{n^3}}} - k\cbraces{1 + 4\frac{k - 2}{k^2}\cdot\frac{4\pi^2}{2n^2} + o\cbraces{\frac{1}{n^3}} + \frac{4}{k^2}\cbraces{\frac{4\pi^2}{2n^2} + o\cbraces{\frac{1}{n^3}}}^2}^{1/2}} \\
	&= \frac{1}{2}\sbraces{k + \frac{4\pi^2}{n^2} + o\cbraces{\frac{1}{n^3}} - k\cbraces{1 + \frac{4\pi^2}{kn^2}\cbraces{1 - \frac{2}{k}} + o\cbraces{\frac{1}{n^3}}}} \\
	&= \frac{1}{2}\sbraces{\frac{4\pi^2}{n^2} - \frac{4\pi^2}{n^2} + \frac{8\pi^2}{kn^2} + o\cbraces{\frac{1}{n^3}}}
	= \frac{4\pi^2}{kn^2} + o\cbraces{\frac{1}{n^3}}
\end{align*}

Therefore, we obtain
\begin{eqnarray*}
	\lambda_{\text{min}}^+(\cG_{RC}) \sim  \Theta \left(  \dfrac{1}{n^2k}  \right).
\end{eqnarray*}